%% file: Simpson2024ParallelizableArXiv.tex
\RequirePackage[OT1]{fontenc}
\documentclass[letterpaper, 10 pt, conference]{ieeeconf}
 
\IEEEoverridecommandlockouts

\usepackage{cite}
\usepackage{comment}
\usepackage{amsmath,amssymb,amsfonts, amsthm}
\usepackage{relsize}
\usepackage{algorithmic}
\usepackage{graphicx}
\usepackage{textcomp}
\usepackage{xcolor}
\usepackage{mathtools}
\usepackage{bm}
\usepackage{optidef}
\usepackage{pgf}
\usepackage{float}
\usepackage{booktabs}
\usepackage{multirow}
\usepackage{tabularx}
\usepackage{dutchcal}
\DeclareMathAlphabet{\matholdcal}{OMS}{cmsy}{m}{n}

\newcommand{\norm}[1]{\left\lVert#1 \right\rVert}
\newcommand{\wnorm}[2]{\norm{#1}_{#2}}

\newcommand{\Model}{
}
\newcommand{\rr}{\mathbb{R}}

\newcommand{\ee}{\mathbb{E}}
\newcommand{\eeo}[1]{\underset{#1}{\ee}}

\newcommand{\tr}{\mathrm{Tr}}
\newcommand{\MHE}{\textrm{MHE}}
\newcommand{\trnorm}[1]{{\left\vert\kern-0.25ex\left\vert\kern-0.25ex\left\vert #1 
		\right\vert\kern-0.25ex\right\vert\kern-0.25ex\right\vert}}
\newcommand{\Nn}{\matholdcal{N}}

\newcommand{\Ll}{\matholdcal{L}}
\newcommand{\Gg}{\matholdcal{G}}

\newcommand{\Uu}{\matholdcal{U}}
\newcommand{\Aa}{\matholdcal{A}}

\newcommand{\Pp}{\matholdcal{P}}
\newcommand{\meanx}{\bar{s}}
\newcommand{\dirfigures}{Figures}

\def\BibTeX{{\rm B\kern-.05em{\sc i\kern-.025em b}\kern-.08em
		T\kern-.1667em\lower.7ex\hbox{E}\kern-.125emX}}

\newtheorem{thm}{Theorem}
\newtheorem{assum}{Assumption}

\newtheorem{remark}{Remark}

\begin{document}
	
	\title{\LARGE \bf
		Parallelizable Parametric Nonlinear System Identification via tuning of a Moving Horizon State Estimator}
	\author{{
			L\'eo Simpson$^{1, 2}$,
			Jonas Asprion$^1$,
			Simon Muntwiler$^3$,
			Johannes K\"{o}hler$^3$,
			Moritz Diehl$^{2,4}$
		}
		\thanks{$^1$ Research and Development, Tool-Temp AG, Switzerland,
			\tt{\small leo.simpson@tool-temp.ch. }}
		\thanks{$^2$ Department of Mathematics, University of Freiburg, 79104 Freiburg, Germany}
		\thanks{$^3$  Institute for Dynamic Systems and Control, ETH Zürich,
			8092 Zürich, Switzerland,
			\tt{\small\{simonmu, jkoehle\}@ethz.ch}}
		\thanks{
			$^4$ Department of Microsystems Engineering (IMTEK),
			University of Freiburg, 79110 Freiburg, Germany
			\tt{\small moritz.diehl@imtek.uni-freiburg.de}}
		\thanks{This research was supported by the EU via ELO-X 953348 and by the Bosch Research Foundation im Stifterverband.}}
	\maketitle
	\begin{abstract}
		\input{sections/abstract}
	\end{abstract}
	
	\begin{keywords}
		Estimation,
            Nonlinear Systems Identification,
            Predictive control for nonlinear systems
	\end{keywords}
	
	 \section{Introduction}
	 \input{sections/intro}

	 \section{Problem Statement}\label{section_problem}
	 \input{sections/setup}
	
	 \section{Prediction Error Method with Moving Horizon State Estimation}\label{section_method}
	 \input{sections/method}

	\section{Special case: Linear Time-Invariant Systems}\label{section_consistencySmall}
	\input{sections/consistencySmall}

	 \section{Numerical Examples}\label{section_numerical}
	 \input{sections/numerical}

	 \section{Conclusion and Outlook}\label{section_conclusion}
	 \input{sections/conclusion}

	\bibliographystyle{IEEEtran}
	\bibliography{biblio}
	
\end{document}

%% file: sections/abstract.tex
This paper introduces a novel
optimization-based approach for parametric nonlinear system identification.
Building upon the prediction error method framework, traditionally used for linear system identification,
we extend its capabilities to nonlinear systems.
The predictions are computed using a moving horizon state estimator with a constant arrival cost.
Eventually, both the system parameters and the arrival cost are estimated by minimizing the sum of the squared prediction errors.
Since the predictions are induced by
the state estimator, the method can be viewed as the tuning of a state estimator, based on its predictive capacities.
The present extension of the prediction error method not only enhances performance for nonlinear systems but also
enables learning from multiple trajectories
with unknown initial states,
broadening its applicability in practical scenarios.
Additionally, 
the novel formulation leaves room for the design of efficient and parallelizable optimization algorithms,
since each output prediction only depends
on a fixed window of past actions and measurements.
In the special case of linear time-invariant systems,
we show an important property of the proposed method which suggests asymptotic consistency under reasonable assumptions.
Numerical examples illustrate the effectiveness and practicality of the approach, and
one of the examples also highlights the necessity for the arrival cost.

%% file: sections/intro.tex
Accurate predictive models are indispensable for reliably forecasting system behaviors,
notably in applications involving the deployment of Model Predictive Control (MPC)~\cite{Rawlings2017}.
Essential to this task is the identification of models using both measurement data and physical knowledge about the system.
Additionally,
beyond having a good predictive model,
the implementation of a reliable state estimator is crucial.
In this paper, we present a method that addresses both requirements simultaneously for nonlinear systems.

In the realm of system identification,
conventional methods such as~\cite{Verhaegen1994, Van1994} often lack the capability to
enforce specific structural constraints on the model.
Parametric system identification overcomes this limitation~\cite{Ljung1999},
with the primary tool being the Prediction Error Method (PEM)~\cite{Ljung2002}.
For state-space models, the PEM uses a state estimator to make output predictions, and the errors of these predictions are minimized over the model parameters.
Therefore, the PEM connects tuning a state estimator and identifying a model.
As it is shown in~\cite{Simpson2023}, this method also offers the possibility of jointly estimating the model parameters and noise covariance matrices, which are key to the design of state estimators.
Traditionally, the PEM is designed for linear system identification, using a Kalman Filter (KF)~\cite{Kalman1960} for the state estimator.
However, the design of state estimators for nonlinear systems is non-trivial and the standard methods are
the Extended Kalman Filter (EKF)~\cite{McGee1985} and the Moving Horizon Estimation (MHE) method~\cite[Chapter 4]{Rawlings2017}.
Recent research has extensively analyzed MHE's capabilities~\cite{alessandri2008moving, Baumgaertner2019, schiller2023lyapunov},
but it has not used it within the PEM framework yet.

The main drawback of the PEM is that it involves solving a possibly complex nonlinear optimization problem even when the system is linear~\cite{Ljung2002}, and the complexity scales with the data-length.
To apply a gradient-based approach to this optimization problem, the derivatives of the state estimates should be computed.
Because of the recursive aspect of the KF or EKF,
this typically involves the propagation of derivatives along the whole input-output dataset,
and these expensive computations can not be made in parallel.
Therefore, a different choice of state-estimator might not only extend the PEM to nonlinear settings,
but also leave room for more efficient optimization algorithms.

\subsubsection*{Contribution}
In this paper, we extend the classical PEM framework to nonlinear systems by
substituting the linear observer for an MHE scheme with a constant arrival cost.
In addition to the typically enhanced performance of MHE over EKF~\cite[Chapter 4]{Rawlings2017},
another advantage stems from its non-recursive nature when the arrival cost is chosen to be constant, enabling state estimation based solely on a finite amount of past inputs and outputs.
Moreover, because of this feature, the proposed approach facilitates the handling of both single and multiple trajectory data.
As mentioned in~\cite{zheng2020},
considering multiple trajectories has the benefit of enabling the identification of unstable systems.
While this method can be viewed as tuning an MHE-based state estimator,
it also has the capability of recovering the true parameters of a model.
Some theoretical guarantees are provided regarding this aspect
in the special case of linear time-invariant systems.
Finally, we provide two numerical examples to showcase the performance
of the proposed method.
The first example also shows by counter-example that
the constant arrival cost is necessary since removing it results in a method that
may fail to recover the true system parameters.

\subsubsection*{Related Work}
Our work builds upon the existing PEM framework~\cite{Ljung2002} but extends its applicability to nonlinear systems and addresses the limitations associated with single but long trajectory data.
In~\cite{Valluru2017}, a similar extension is introduced using the EKF for state estimation.
By weighting the prediction errors with the predicted output covariance,
this method is equivalent to Maximum Likelihood Estimation (MLE)~\cite{astrom1979maximum}
in the special case of linear systems.
In that latter case, the theoretical guarantees are similar to the ones in~\cite{Ljung2002}.
In particular, they only hold in the case of learning from an infinitely long single trajectory.
On the other hand, the treatment of multi-trajectory system identification has been explored in several studies~\cite{zheng2020, tu2022learning, dean2020sample},
but on the restricted case of independent trajectories of a linear system.
In~\cite{dean2020sample}, the system identification task is turned into a regression problem
with independent samples.
In~\cite{tu2022learning},
single-trajectory and multiple-trajectory learning are treated in a unified way.
Notably, the authors draw the conclusion that learning from many small trajectories is more efficient
than learning from a single long trajectory.

Finally,
in~\cite{Muntwiler2022} the authors tune a linear MHE scheme using the sensitivities of the MHE solution with respect to system parameters.
Other gradient-based approaches have been proposed for tuning some parameters entering the MHE cost
~\cite{esfahani2021reinforcement, esfahani2023learning}.
Since the purpose was not system identification,
the authors do not treat the question of recovering the parameters of a model.
Furthermore, since a dynamic arrival cost is used, parallelization of the optimization algorithm is not possible.

\subsubsection*{Outline}
The remainder of this paper is organized as follows:
Section~\ref{section_problem}
defines the precise settings of the identification task and
Section~\ref{section_method} presents the proposed solution.
In Section~\ref{section_consistencySmall},
an important property is proven for the special case of Linear Time-Invariant (LTI) Systems.
In Section~\ref{section_numerical},
two simulation examples are provided,
and the performance of the proposed method is illustrated and compared to
a more naive approach.
Lastly, some conclusions are drawn and future research directions
are suggested in Section~\ref{section_conclusion}.

%% file: sections/setup.tex
The task is to identify a model for the evolution of a system, from which measurement data is available.
We assume that a parametric model is known, and
the parameter $\theta \in \rr^{n_{\theta}}$ corresponding to the system is to be estimated from measurement data.
Such a parametric model can, for example, come from physics-based knowledge about the system.
In this paper, the parametric model takes the following nonlinear state-space form,
with additive noise:
\begin{align}\label{equation_model}
    \begin{split}
        x_{t+1} &= f_{\Model}(x_t, u_t; \theta) + w_t, \\
        y_t &= g_{\Model}(x_t; \theta) + v_t,
    \end{split}
\end{align}
where $x_t \in \rr^n$ represents the state of the system,
$u_t \in \rr^q$ are the inputs,
and $y_t \in \rr^p$ are the outputs.
The unknown parameter $\theta$ is assumed to be in a set $\Pp$ defined as follows
\begin{equation}
	\Pp \coloneqq \left\{ \theta \in \rr^{n_\theta} \big\vert  h_{\Pp}(\theta) \leq 0  \right\},
\end{equation}
for some continuous function $h_{\Pp}(\cdot)$.
The process and measurements noise $w_t$ and $v_t$
are assumed to be independent and identically distributed, with zero mean.
Their covariance matrices $Q(\theta)$ and $R(\theta)$
may also depend on the unknown parameter $\theta$.
For the method to be well defined and to ensure that derivative-based approaches can be applied for our optimization-based method, 
we assume that $Q(\theta)$ and $R(\theta)$ are positive definite for any $\theta \in \Pp$, and that the functions $f$, $g$, $Q$ and $R$
are twice continuously differentiable with respect to all their arguments.

The available data is a set of different input and output sequences collected at different times, from measurements of the system:
\begin{align}
    \begin{split}
        &u^{(i)} \coloneqq  \begin{bmatrix} 
            u_{t^{(i)}} \\ \dots \\  u_{t^{(i)}+m-1} \end{bmatrix}, \quad
        y^{(i)} \coloneqq  \begin{bmatrix} 
            y_{t^{(i)}+1} \\ \dots \\  y_{t^{(i)}+m} \end{bmatrix},
    \end{split}
\end{align}
where for $i=1, \dots, N$, the $i$-th data sequence is collected from time $t^{(i)}$ and is of length $m$.
To keep the notation simple, we assume that $m$ is independent of $i$. 
Note that since we do not assume that the sequences are independent, the corresponding time windows might intersect.

\begin{remark}\label{remark:oneTrajectory}
	The case of a single trajectory $u_0, y_1, \dots, u_{\tilde{N}-1}, y_{\tilde{N}}$
    is covered by $t^{(i)} = i$ and $N = \tilde{N}-m$.
\end{remark}

%% file: sections/method.tex
\begin{figure}
	\vspace{0.3cm}
	\begin{center}
		\includegraphics[width=0.95\linewidth]{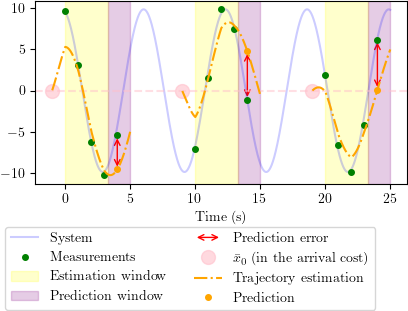}
	\end{center}
	\vspace{-0.5cm}
	\caption{
	Trajectory estimation (orange dashed line) using MHE with a constant arrival cost and a horizon of $4$ steps and one-step-ahead prediction errors (red arrows).
	The system is a harmonic oscillator $\ddot{x} = - x$ with discrete measurements (green dots).
    The MHE is computed for a model mismatch $\ddot{x} = - \omega^2 x$ with $\omega = 0.7$.
	}\label{figure_method}
 \vspace{-0.5cm}
\end{figure}

The estimation method that we present falls into the PEM framework~\cite{Ljung2002},
where predictions of the last outputs are compared with the corresponding measurement
to construct the prediction errors.
For a better presentation of the method, we use the following shorthand notations:
\begin{align}
	\begin{split}
		&
		z^{(i)} \coloneqq
		\begin{bmatrix}
			y_{t^{(i)}+1} \\ \dots \\  y_{t^{(i)}+m-1}
		\end{bmatrix},
		\quad
		\tilde{y}^{(i)} \coloneqq y_{t^{(i)}+m}.
	\end{split}
\end{align}
For $i=1, \dots, N$, the data $u^{(i)}, z^{(i)}$ is used to estimate the state at time $t^{(i)}+m$.
Ultimately, these state estimates will induce output predictions for time $t^{(i)}+m$,
which are compared to the actual outputs $\tilde{y}^{(i)}$ to construct the prediction errors.
These prediction errors are minimized over the parameters of the model.
The corresponding minimizer gives the desired estimate of the parameter $\theta$.
The approach is illustrated for $m=5$ in Fig.~\ref{figure_method}.

\begin{remark}
To keep the method simple, we choose the squared $L_2$ norm of the prediction errors for the objective function.
This implies, for example, that the matrices $Q$ and $R$ can only be estimated up to a multiplicative constant.
Some variants exist in the literature~\cite{Ljung1999}, notably the important special case of the MLE method~\cite{astrom1979maximum},
where the error is weighted by the information matrix of the output predictions.
\end{remark}

\subsection{Moving Horizon State Estimation}\label{subsection_MHE}

Regarding the state estimator, we use a standard MHE method for nonlinear systems~\cite[Chapter 4]{Rawlings2017}.
In the unconstrained case, the MHE method
consists in minimizing a cost $\Ll_{\MHE}$ over
a state trajectory $\hat{x} \coloneqq \begin{bmatrix} \hat{x}_0 & \cdots & \hat{x}_m \end{bmatrix}^\top$, given some input and output sequences $u$ and $z$.
This cost is defined as follows:
\begin{equation}\label{equation_MHE}
	\begin{aligned}
		&\Ll_{\MHE}(\hat{x}, z, u; \theta, \meanx, \Sigma) \coloneqq \\
		& \qquad
		\wnorm{\hat{x}_0 - \meanx }{{\Sigma}^{-1}}^2
		+ \sum_{k=1}^{m} \wnorm{ z_{k} - g_{\Model}(\hat{x}_k; \theta)}{{R(\theta)}^{-1}}^2\\
		& \qquad
		\phantom{\lVert x_0 - \meanx \rVert_{{\Sigma}^{-1}}^2}
		+ \sum_{k=1}^{m} \wnorm{\hat{x}_{k} - f_{\Model}(\hat{x}_{k-1}, u_k; \theta)}{{Q(\theta)}^{-1}}^2,
	\end{aligned}
\end{equation}
where $\meanx$ and $\Sigma$ parameterize the arrival cost,
and $\norm{x}_M$ denotes the weighted norm $\sqrt{x^\top M x}$
for any positive-definite matrix $M \in \rr^{n\times n}$ and vector $x \in \rr^n$.
Note that one can interpret~\eqref{equation_MHE}  as the maximum a posteriori state trajectory problem,
when a Gaussian prior is assumed on the initial state and noise~\cite{rao2002constrained}.
The solution $\hat{x}_{m}$ of the optimization problem~\eqref{equation_MHE}
is the state estimate at time step $m$ which is used to build a prediction of the next output $\tilde{y}$.

\subsection{Constant Arrival Cost}\label{subsection_arrivalCost}
In the MHE formulation~\eqref{equation_MHE}, the arrival cost corresponds to
a prior knowledge on the states $x_{t^{(i)}}$, and is parameterized by a mean $\meanx$ and a covariance
$\Sigma$.
In the present method, we parameterize $\meanx$ and $\Sigma$
with an additional parameter $\eta \in \Aa$ (for some feasible set $\Aa$) that we optimize jointly with $\theta$:
\begin{align}
	\begin{split}
		\meanx &= \meanx(\eta), \quad
		\Sigma = \Sigma(\eta).
	\end{split}
\end{align}
For the MHE scheme~\eqref{equation_MHE} to be well defined,
we choose a parameterization and a feasible set $\Aa$ such that for any $\eta \in \Aa$, the matrix
$\Sigma(\eta)$ is positive definite.
When a full parameterization of $\meanx, \Sigma$ is desired,
a typical choice is to parameterize
the Cholesky decomposition of $\Sigma^{-1}$.
More precisely, one can choose
$\eta \coloneqq (\meanx, L)$ where $L$ is lower triangular, and
$\Sigma(\eta) \coloneqq {( L L^\top)}^{-1}$.

\begin{remark}
	The estimates of $\eta$ are highly influenced by the input sequence of the available data,
	so the corresponding arrival cost may not be suitable for implementing an MHE scheme on the system when the inputs are distributed differently.
    Instead, a recursive arrival cost is more suited for online state estimation~\cite{alessandri2008moving}.
\end{remark}
Note that a tempting approach would have been to remove this arrival cost since we do not have information about the initial state.
In Section~\ref{section_numerical}, we show by a counter-example that this may result in a biased estimate of $\theta$.

\subsection{Prediction Error Method}
As explained above,
the estimates given by the present method are obtained by
minimizing the output prediction errors.
Since these output predictions are based on an optimization-based state estimate,
the present method can be formulated as a bilevel optimization problem as follows:
\begin{equation}\label{equation_PredErr}
\begin{aligned}
	& 
	\underset{\bm{\hat{x}}, \; 
		\theta \in \Pp, \;
        \eta \in \Aa
	}{\textup{minimize}}
	\quad
	\frac{1}{N} \sum_{i=1}^N \lVert
	\tilde{y}^{(i)} - g_{\Model}(\hat{x}_{m}^{(i)} ; \theta) \rVert^2
	\\
	& \textup{subject to }\\
	&\quad \hat{x}^{(i)}  \in
		\arg \min_{\hat{x}} \Ll_{\MHE}\big(\hat{x}, z^{(i)}, u^{(i)}; \meanx(\eta),  \Sigma(\eta)\big),
	\\
	& \hspace{2cm}  \textup{for} \; i = 1, \dots, N,
\end{aligned}
\end{equation}
where $\bm{\hat{x}} \coloneqq \begin{bmatrix}
    \hat{x}^{(1)} & \cdots & \hat{x}^{(N)}
\end{bmatrix}$
are optimization variables implicitly defined by the argmin constraint.
Note that this formulation differs from the standard PEM since
the predictions only depend on a fixed number of past inputs and outputs.
This aspect makes consistency guarantees more difficult to prove
since these proofs often need the effect of the unknown initial state to dissipate over a long time horizon.
This is why the output predictions are usually computed using
all of the past data.
On the other hand, it may grant computational advantages since the constraints are now independent of each other and sparsity of the optimization problem is increased.

To solve the bilevel optimization problem \eqref{equation_PredErr}, we relax the global optimality condition of the inner optimization problems by their first-order optimality conditions:
\begin{align}\label{equation_PredErrKKT}
	& 
	\underset{\bm{\hat{x}}, \; 
		\theta \in \Pp, \;
        \eta \in \Aa
	}{\textup{minimize}}
	\quad
	\frac{1}{N} \sum_{i=1}^N \lVert
	\tilde{y}^{(i)} - g_{\Model}(\hat{x}_{m}^{(i)} ; \theta) \rVert^2
	\\ \nonumber
	& \textup{subject to }\\ \nonumber
	&\quad \nabla_{\hat{x}} \Ll_{\MHE}\big(\hat{x}^{(i)}, z^{(i)}, u^{(i)}; \meanx(\eta),  \Sigma(\eta)\big)
	= 0
	\\ \nonumber
	& \hspace{2cm}  \textup{for} \; i = 1, \dots, N.
\end{align}
\begin{remark}
	This relaxation is exact if the inner optimization problems are convex.
\end{remark}

The constraints of problems~\eqref{equation_PredErr} and~\eqref{equation_PredErrKKT} being independent, most of the
computational effort required for gradient-based optimization may be parallelized.
More specifically,
the solutions of the inner optimization problems in~\eqref{equation_PredErr}
together with their sensitivities with respect to $(\theta, \eta)$
can be computed in parallel.
Additionally, these inner optimization problems have sparsity patterns similar to the ones encountered in Optimal Control Problems (OCP), thus, they can be solved with a linear complexity in the MHE horizon $m$.
For these reasons, we expect that problems~\eqref{equation_PredErr} and~\eqref{equation_PredErrKKT} can be efficiently solved.
However, investigating new numerical optimization algorithms is outside the scope of this paper and is left for future research.

%% file: sections/consistencySmall.tex
In this section, we prove an important statistical property regarding the capacity of the
presented method to recover the true parameter $\theta^\star$ in the special case of
LTI systems.
The model reads:
\begin{equation}\label{equation_LTIsmall}
    \begin{aligned}
        x_{t+1} &= A(\theta) x_t + B(\theta) u_t + w_t, \\
        y_t &= C(\theta) x_t + v_t,
    \end{aligned}
\end{equation}
where $A(\theta)$, $B(\theta)$ and $C(\theta)$
are matrices of appropriate dimension that depend on the unknown parameter $\theta$.
We also assume the process and measurement noise to be normally distributed:
$w_t \sim \Nn \big(0, Q(\theta) \big)$ and
$v_t \sim \Nn \big(0, R(\theta) \big)$.

As in the paper, the available data is a list of subsequences of inputs and outputs starting from $t=t^{(i)}$.
\begin{assum}\label{assum:iid}
    We assume that the input subsequences $u^{(i)}$ are random variables with identical mean denoted by
$\bar{u}$, and independent from the initial states $x_{t^{(i)}}$ and from the process and measurement noise.
\end{assum}
Note that Assumption~\ref{assum:iid} holds for example in the case of independent and identically distributed inputs $u_{t}$ or in the case of independent trajectories.

As discussed in Section~\ref{subsection_arrivalCost} the arrival cost is fully parameterized.
We denote
the output predictions based on
the inputs $u^{(i)}$, the outputs $z^{(i)}$ and
the parameters $\left(\theta, \meanx, \Sigma \right)$, by
$\hat{y}(u^{(i)}, z^{(i)}; \theta, \meanx, \Sigma)$.
These are given by
the state estimate computed as discussed in Section~\ref{subsection_MHE},
multiplied by the output matrix $C_{\theta}$.
We formulate the optimization problem~\eqref{equation_PredErr}
in a condensed way:
\begin{equation}\label{equation_condensed}
        \underset{\substack{
                \theta \in \Pp, \meanx, \Sigma \succeq 0
        }}{\textrm{minimize}}
        \,
        \sum_{i=1}^N
        \norm{
            \tilde{y}^{(i)} \! \! -
            \hat{y}\! \left(u^{(i)}, z^{(i)}; \theta, \meanx, \Sigma \right)\!
        }^2
        \! \! \! \eqqcolon \! V_N(\bm{u}, \bm{y}; \theta, \eta)
\end{equation}
\begin{thm}
    Assume the data is generated by the LTI system~\eqref{equation_LTIsmall} under Assumption~\ref{assum:iid}
    and with parameters $\theta = \theta^\star$.
    Then there exists some $\eta^\star_N$
    such that $(\theta^\star, \eta^\star_N)$
    minimizes the expected value of
    the objective function in~\eqref{equation_condensed},
    i.e.
    \begin{equation}\label{equation_consistencySmall}
        (\theta^\star, \eta_N^\star) \in \arg\min_{\theta, \eta} \;  \underset{}{\eeo{\bm{u}, \bm{y}}} \bigg[ V_N(\bm{u}, \bm{y}, \theta, \eta) \bigg].
    \end{equation}
\end{thm}
\begin{remark}
    This result does not directly imply the asymptotic consistency of the presented method,
    i.e.\ the convergence, as the number of data sequences $N$ goes to infinity,
    of the parameter estimate to the parameter that was used to generate the data.
    However, it provides positive insights regarding the ability to recover
    the true parameters of the system.
    Indeed, equation~\eqref{equation_consistencySmall} suggests that the parameter estimate is close to the ground truth $\theta^\star$
    if the objective function gets close to its expected value,
    and if the latter has a unique minimizer.
    A comprehensive proof of asymptotic consistency with appropriate assumptions
    is left for future research.
\end{remark}

\begin{proof}
    \item
    \subsubsection{Notation}
    For compactness, we define the following vectors:
    \begin{equation*}
    \begin{aligned}
        s^{(i)} \coloneqq x_{t^{(i)}},
        \quad
        e^{(i)} \coloneqq \begin{bmatrix}w^{(i)} \\ v^{(i)} \end{bmatrix},
    \end{aligned}
    \end{equation*}
    where $w^{(i)}$ and $v^{(i)}$ stack the process and
    measurement noise.
    Note that the initial states $s^{(i)}$ and the noises
    $e^{(i)}$
    are unknown to the estimation routine.

    \subsubsection{Definition of $\eta^\star_N$}
    We define $\meanx_N$ and $\Sigma_N$ as 
    the empirical mean and covariance of the (unknown) initial states:
    \begin{equation*}
        \begin{aligned}
            \meanx_N &\coloneqq \frac{1}{N} \sum_{i=1}^N \eeo{s^{(i)}}\left[ s^{(i)} \right], \\
            \Sigma_N &\coloneqq \frac{1}{N} \sum_{i=1}^N
            \eeo{s^{(i)}}\left[ (s^{(i)} - \meanx_N) {(s^{(i)} - \meanx_N)}^\top\right].
        \end{aligned}
    \end{equation*}
    Since $\meanx$ and $\Sigma$ are assumed to be fully parameterized, we define $\eta_N^\star$ accordingly.
    By construction, the following holds:
    \begin{equation}\label{equation_MeanSigmaProp}
    	\begin{aligned}
    		\sum_{i=1}^N \eeo{d^{(i)}}\left[ d^{(i)} \right] &= 0, \\\quad
    		\sum_{i=1}^N \eeo{d^{(i)}}\left[ d^{(i)}d^{(i) \top} \right]  &= N\Sigma_N,
    	\end{aligned}
    \end{equation}
    with $d^{(i)}\coloneqq s^{(i)} - \meanx_N$.

    \subsubsection{Stochastic interpretation of the predictions}
Since the model is linear and the noise is Gaussian,
	the predictions
	$\hat{y}\left(u^{(i)}, z^{(i)} ; \theta^\star, \meanx_N, \Sigma_N \right)$
	are the KF predictions with a Gaussian prior knowledge of the initial state
	$s^{(i)}$ that we write $\Gg^{(i)} \coloneqq
	\left( s^{(i)} \sim \Nn(\meanx_N, \Sigma_N) \right)$.
	Notably, it implies that the predictions are
	exactly the expected values of $\tilde{y}^{(i)}$
	conditioned by $u^{(i)}, z^{(i)}$ and $\Gg^{(i)}$:
	\begin{equation*}
		\hat{y}(u^{(i)}, z^{(i)}, \theta^\star, \meanx_N, \Sigma_N) = \eeo{e^{(i)}, s^{(i)}}\left[
		\tilde{y}^{(i)} \big\vert u^{(i)}, z^{(i)}, \Gg^{(i)}
		\right],
	\end{equation*}
	which implies:
	\begin{equation*}
		\begin{aligned}
			 \eeo{e^{(i)}, s^{(i)}} \left[ 
			\varepsilon^{(i)} r^{(i) \top} \big\vert u^{(i)}, z^{(i)}, \Gg^{(i)}
			\right] = 0,
		\end{aligned}
	\end{equation*}
	where we defined:
	\begin{equation*}
		r^{(i)} \coloneqq \begin{bmatrix}1 \\ u^{(i)} \\ z^{(i)} \end{bmatrix}, \quad
		\varepsilon^{(i)} \coloneqq
		\tilde{y}^{(i)} - \hat{y}(u^{(i)}, z^{(i)};\theta^\star, \meanx_N, \Sigma_N).
	\end{equation*}
Furthermore, this implies, using the law of total expectation:
	\begin{equation*}
		\begin{aligned}
			& \eeo{e^{(i)}, s^{(i)}} \left[ 
			\varepsilon^{(i)} r^{(i) \top}
			\big\vert u^{(i)}, \Gg^{(i)} \right] \\
			& \hspace{1cm} =\eeo{z^{(i)}} \left[ 
			\eeo{e^{(i)}, s^{(i)}} \left[ 
			\varepsilon^{(i)} r^{(i) \top} \big\vert u^{(i)}, z^{(i)}, \Gg^{(i)}
			\right] \big\vert  u^{(i)}, \Gg^{(i)} \right] \\
			&\hspace{1cm} = \eeo{z^{(i)}} \left[ 
			0 \big\vert  u^{(i)}, \Gg^{(i)} \right] = 0.
		\end{aligned}
	\end{equation*}
Note that one can also define the Gaussian prior using $d^{(i)}$:
$\Gg^{(i)} \coloneqq
\left( d^{(i)} \sim \Nn(0, \Sigma_N) \right)$, and the latter equation can be manipulated using the total expectation law as follows:
\begin{equation}\label{equation_conditionalcorrelations}
		\begin{aligned}
		&\eeo{d^{(i)}}  \left[ 
			\eeo{e^{(i)}}  \left[ 
			\varepsilon^{(i)} r^{(i) \top} \big\vert u^{(i)}, d^{(i)}
			\right] \big\vert u^{(i)}, \Gg^{(i)} \right]
   \\
    & \hspace{1cm} =
			\eeo{e^{(i)}, d^{(i)}} \left[ 
			\varepsilon^{(i)} r^{(i) \top}
			\big\vert u^{(i)}, \Gg^{(i)} \right] = 0.
		\end{aligned}
\end{equation}

    \subsubsection{Covariance between prediction errors and past data}
    The outputs $z^{(i)}$ and $\tilde{y}^{(i)}$
    are linear functions of $s^{(i)}$, $e^{(i)}$, and $u^{(i)}$.
    Since $\varepsilon^{(i)}$ is linear in
    $z^{(i)}$, $u^{(i)}$, $\meanx_N$ and $\tilde{y}^{(i)}$,
    it is also a linear function of
    $\meanx_N$, $s^{(i)}$, $e^{(i)}$ and $u^{(i)}$.
    Finally, it can also be viewed as an affine function of
    $d^{(i)}$, $e^{(i)}$ and $u^{(i)}$.
    Therefore, there exist some matrices $M_1$ and $M_2$
    such that the following holds for $i=1, \dots, N$:
    \begin{equation*}
        \begin{aligned}
            \varepsilon^{(i)} r^{(i) \top} =
            M_1
            \begin{bmatrix}
                d^{(i)} \\
                e^{(i)} \\
                u^{(i)} \\
                1
            \end{bmatrix}
            \begin{bmatrix}
                d^{(i) \top} &
                e^{(i) \top} &
                u^{(i) \top} &
                1
            \end{bmatrix}
            M_2.
        \end{aligned}
    \end{equation*}
Thus, using the independence between
    $e^{(i)}$ and $( u^{(i)}, d^{(i)})$:
\begin{equation*}
        \begin{aligned}
        \eeo{e^{(i)}} \left[ 
            \varepsilon^{(i)} r^{(i) \top} \big\vert u^{(i)}, d^{(i)}
            \right]
        &=
        \varphi_u \left(u^{(i)}  \right)
        + \varphi_{uu} \left(u^{(i)} u^{(i) \top}  \right)
        \\ &\hspace{-3cm}
         +\varphi_d \left(d^{(i)}  \right)
        + \varphi_{du} \left(d^{(i)} u^{(i) \top}  \right)
         + \varphi_{dd} \left(d^{(i)} d^{(i) \top}  \right) + \varphi_{0},
        \end{aligned}
    \end{equation*}
    for some linear matrix functions $\varphi_u$, $\varphi_{uu}$, $\varphi_d$, $\varphi_{du}$, $\varphi_{dd}$ and some matrix $\varphi_{0}$.
    Inserting this into equation~\eqref{equation_conditionalcorrelations}:
\begin{equation*}
        \begin{aligned}
        \varphi_u \left(u^{(i)}  \right)
        + \varphi_{uu} \left(u^{(i)} u^{(i) \top}  \right)
         + \varphi_{dd} \left(\Sigma_N  \right) + \varphi_{0}
         = 0,
        \end{aligned}
    \end{equation*}
    which holds for all $u^{(i)}$.
    This implies that $\varphi_u$ and $\varphi_{uu}$ are null and $\varphi_0 = - \varphi_{dd} \left( \Sigma_N \right)$.
    Overall, the following holds:
    \begin{equation*}
        \begin{aligned}
        \eeo{e^{(i)}} \!\left[ 
            \varepsilon^{(i)} r^{(i) \top} \big\vert u^{(i)}, d^{(i)}
            \right]
        &=
        \varphi_d \left(d^{(i)}  \right)
        + \varphi_{du} \left(d^{(i)} u^{(i) \top}  \right) 
        \\& \hspace{0.5cm}
         + \varphi_{dd} \left(\!d^{(i)} d^{(i) \top} - \Sigma_N \right).
        \end{aligned}
    \end{equation*}
    After defining
    $\tilde{\varphi}_d(d^{(i)}) \coloneqq \varphi_d \left(d^{(i)}  \right)  + \varphi_{du} \left(d^{(i)} \, \bar{u}^\top \right)$, where $\bar{u}$ is according to Assumption~\ref{assum:iid}, the following holds:
    \begin{equation*}
        \begin{aligned}
        \eeo{e^{(i)}, u^{(i)}} \!\left[ 
            \varepsilon^{(i)} r^{(i) \top} \big\vert d^{(i)}
            \right]
        &=
         \tilde{\varphi}_d \left(d^{(i)}  \right) 
         + \varphi_{dd} \left(\!d^{(i)} d^{(i) \top}\! \! - \! \Sigma_N \right).
        \end{aligned}
    \end{equation*}
    Finally, combining the last equation and~\eqref{equation_MeanSigmaProp}:
    \begin{align}\label{equation_AlmostFinalResult}
        \eeo{\bm{u}, \bm{y}} \left[  \sum_{i=1}^N 
            \varepsilon^{(i)} r^{(i) \top} \right]
        &= 
        \! \eeo{d^{(1)}\!\!, \dots, d^{(N)}} \!
        \left[
        \sum_{i=1}^N
        \eeo{e^{(i)}, u^{(i)}} \left[ 
            \varepsilon^{(i)} r^{(i) \top} \vert d^{(i)}
            \right] \right] \nonumber \\ \nonumber
        &\hspace{-2cm}=
        \eeo{\bm{d}} \left[ \tilde{\varphi}_d \left(\sum_{i=1}^N d^{(i)}  \right) 
         + \varphi_{dd} \left(\sum_{i=1}^N d^{(i)} d^{(i) \top}-N\Sigma_N  \right) \right] \\
         \nonumber
         &\hspace{-2cm}=
        \tilde{\varphi}_d \left(\sum_{i=1}^N \!\eeo{d^{(i)}}\!\!\left[ d^{(i)}\!\right] \! \right) 
         \!+ \! \varphi_{dd} \! \left(
            \sum_{i=1}^N\! \eeo{d^{(i)}}\!\!\left[ d^{(i)}  d^{(i) \top}\!\right]
            \! \! -\! \! N\Sigma_N 
             \right)  \\  
         &\hspace{-2cm}=
        \tilde{\varphi}_d \left(0  \right) 
         + \varphi_{dd} \left(0 \right)
        = 0. 
    \end{align}

    \subsubsection{Linearity of the Prediction Errors}
    For any $\theta, \eta$, the predictions
    $\hat{y}\left(u^{(i)}, z^{(i)}; \theta, \meanx(\eta), \Sigma(\eta) \right) $
    are affine in $u^{(i)}$ and $z^{(i)}$, hence,
    these are linear in $r^{(i)}$ (but nonlinear in $(\theta, \eta)$).
    Therefore, for any $\theta, \eta$, 
    there exists a matrix $\Phi(\theta, \eta)$
    such that for $i=1, \dots, N$:
    \begin{equation*}
        \hat{y}\left(u^{(i)}, z^{(i)}; \theta, \meanx(\eta), \Sigma(\eta) \right) 
        = \Phi(\theta, \eta) r^{(i)}.
    \end{equation*}
    Hence:
    \begin{equation*}
        \begin{aligned}
      &\tilde{y}^{(i)} - \hat{y}\left(u^{(i)}, z^{(i)} ; \theta, \meanx(\eta), \Sigma(\eta) \right) \\
      &  \hspace{1cm} = \tilde{y}^{(i)} - 
        \hat{y}\left(u^{(i)}, z^{(i)};\theta^\star, \meanx_N, \Sigma_N\right) +\Delta(\theta, \eta) r^{(i)}\\
        & \hspace{1cm} = \varepsilon^{(i)}
        +\Delta(\theta, \eta) r^{(i)},
        \end{aligned}
    \end{equation*}
    where $\Delta(\theta, \eta) \coloneqq
    \Phi(\theta, \eta) - \Phi(\theta^\star, \eta^\star_N)$.
    \subsubsection{Conclusion}
    Finally, using the last equation and~\eqref{equation_AlmostFinalResult}:
    \begin{equation*}
        \begin{aligned}
            &\eeo{\bm{u}, \bm{y}}\big[  V_N(\bm{u}, \bm{y}; \theta, \eta) \big] \\
            &= \sum_{i=1}^N
            \eeo{u^{(i)}, y^{(i)}} \! \left[ 
                \norm{
                \tilde{y}^{(i)} - \hat{y}\left(u^{(i)}, z^{(i)} ; \theta, \meanx(\eta), \Sigma(\eta) \right)
                }^2
             \right], \\
             &=
            \sum_{i=1}^N
            \eeo{u^{(i)}, y^{(i)}} \! \left[ 
                \norm{
                    \varepsilon^{(i)} +
                    \Delta(\theta, \eta) r^{(i)}
                }^2
             \right],\\
             &\geq 
            \sum_{i=1}^N
            \eeo{u^{(i)}, y^{(i)}} \! \bigg[ 
                \norm{\varepsilon^{(i)}}^2
              +  2 \tr\left( \varepsilon^{(i)} r^{(i)\top} {\Delta(\theta, \eta)}^\top \right)
             \bigg],\\
             &=
             \eeo{\bm{u}, \bm{y}}\left[ V_N(\bm{u}, \bm{y}; \theta^\star, \eta^\star_N) \right]
             \!+\!2 \tr \!\left(\!\!\Delta(\theta, \eta)\!^\top\!\!
            \eeo{\bm{u}, \bm{y}} \! \!\left[ 
                \sum_{i=1}^N \varepsilon^{(i)} r^{(i)\top}\!
             \right] \!
             \right)
             \! \!, \\
            &= 
            \eeo{\bm{u}, \bm{y}}\left[ V_N(\bm{u}, \bm{y}; \theta^\star, \eta^\star_N) \right]
            .
        \end{aligned}
    \end{equation*}
Hence,
$\eeo{\bm{u}, \bm{y}}\left[ V_N(\bm{u}, \bm{y}; \theta^\star, \eta^\star_N) \right] \leq \eeo{\bm{u}, \bm{y}}\big[ V_N(\bm{u}, \bm{y}; \theta, \eta) \big]$
for any $\theta, \eta$, which concludes the proof.

\end{proof}

%% file: sections/numerical.tex
In this section, we showcase the consistency and practicality of the proposed approach through numerical simulation. 
We study two examples,
where, by applying our method to noisy data,
the model parameters are recovered.
The first example also serves as a counter-example for a naive, yet tempting alternative approach.
Eventually, we consider a more practical example where the system is nonlinear and the noise is non-Gaussian, reflecting real-world scenarios.
Regarding the numerical optimization,
we used the nonlinear optimization solver
IPOPT~\cite{Waechter2006} together with the package
CasADi~\cite{Andersson2019} for the symbolic expressions
and automatic differentiation,
and the sparse linear solver MUMPS\@.

\subsection{A Simple LTI Case Study}
As a first case study, we consider the following single-state LTI system:
\begin{equation}\label{equation_counterex}
	\begin{aligned}
		x_{k+1} &= \theta x_k + w_k,\\
		y_k &= x_k + v_k,
	\end{aligned}
\end{equation}
where $x_k$ is the state, $y_k$ is the output,
$\theta$ is the parameter to be estimated,
and $w_k$, $v_k$ are the process and measurement noise,
which are assumed to be independent and normally distributed
$w_k, v_k \sim \Nn(0, 1)$.
As in Remark~\ref{remark:oneTrajectory}, we choose the setting of a single trajectory,
by choosing $t^{(i)} = i$, and we consider subsequences of the trajectories of size $m$.
Overall, the presented method summarized by the optimization problem~\eqref{equation_PredErr}
takes the following form:
\begin{equation}\label{equation_PredErrExample}
	\begin{aligned}
		& 
		\underset{
			\theta, \lambda, \meanx,
			\bm{\hat{x}}
		}{\textup{minimize}}
		\qquad
		\frac{1}{N} \sum_{i=1}^N \,
		{\left( y_{i+m} - \hat{x}^{(i)}_{m}\right)}^2
		\\
		& \textup{subject to }
		\quad \hat{x}_{0:m}^{(i)}  =
			\arg \min_{x_{0:m}} \,
			\Ll_{\MHE}^{(i)}\big(\hat{x}; \meanx,  \lambda \big),
	\end{aligned}
\end{equation}
where
$\Ll_{\MHE}^{(i)}\big(\hat{x}; \meanx,  \lambda \big) = 
\lambda {\left( \meanx - x_0 \right)}^2 
 + \sum_{k=1}^{m-1} \left( y_{i+k} - x_{k} \right)^2
+ \sum_{k=1}^{m} \left( x_{k} - \theta x_{k-1} \right)^2$,
and
$\lambda$ corresponds to $\Sigma^{-1}$.
In addition, we show by counter-example that while removing the constant arrival cost of the MHE seems tempting,
it results in a biased parameter estimate.
Such a method would correspond to setting
$\lambda=0$ in problem~\eqref{equation_PredErrExample}.

For the numerical experiment, we generate $10$ data instances according to model~\eqref{equation_counterex},
with the parameter $\theta^\star = 0.8$, with different noise realization
and with $x_0 = 0$.
We then run the estimation method~\eqref{equation_PredErrExample} for each data instance,
with $m=3$ and for different data lengths $N = 10^2, \dots, 10^6$.
For comparison, we also run the estimation method with $\lambda=0$ as described above.

\begin{figure}
	\vspace{0.2cm}
	\begin{center}
		\includegraphics[width=0.95\linewidth]{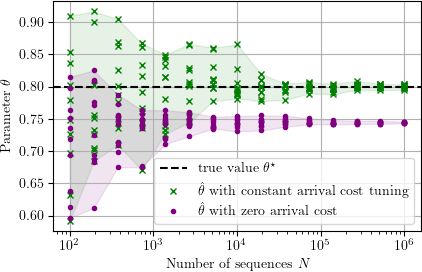}
	\end{center}
	\vspace{-0.5cm}
	\caption{Parameter estimates with and without constant arrival cost tuning for the toy example~\eqref{equation_counterex}.
  For each data instance,
  the estimates given by the proposed approach are represented by a green cross, and the ones given by the naive approach are the purple dots.
  The shaded areas depict, for each data length and for both methods, the interval in which the $10$ parameter estimates lie.}\label{fig_LTIex}
 \vspace{-0.5cm}
\end{figure}

The results are depicted in Fig.~\ref{fig_LTIex}.
They clearly show that for the present example,
the method with zero arrival cost is biased,
while the proposed method is consistent.
This result may seem counter-intuitive at first,
but a deeper analysis would explain the bias.
Since the measurements are noisy and
the states $x_k$ are distributed around zero,
the outputs $y_{i+m}$ are, on average, closer to zero
than their predictions based on measurements only and for $\theta=\theta^\star$.
As a result, shrinking down the state estimates induces lower prediction error on average.
This is achieved for $\theta < \theta^\star$, hence the bias when no arrival cost is considered.
Another interesting remark is that the method with constant arrival cost tuning seems
to require more data to reach its limit.
Since the proposed method additionally estimates $\meanx, \lambda$,
it is not surprising that more data is required.

\subsection{Lorenz Attractor}
\begin{figure}
	\vspace{-0.2cm}
	\begin{center}
		\includegraphics[width=0.95\linewidth]{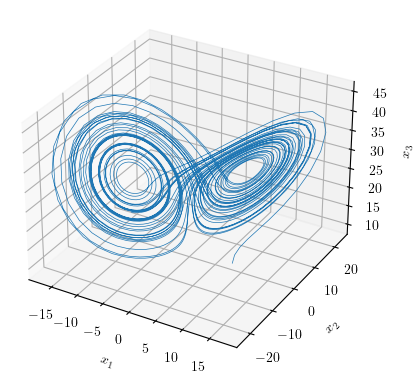}
	\end{center}
	\vspace{-0.5cm}
	\caption{Simulation of the Lorenz attractor for parameters
	$\theta = \begin{bmatrix}10 & 30 \end{bmatrix}^\top$.
	}\label{fig_Lorenz}
 \vspace{-0.5cm}
\end{figure}

\begin{figure}
	\vspace{0.3cm}
	\begin{center}
		\includegraphics[width=0.95\linewidth]{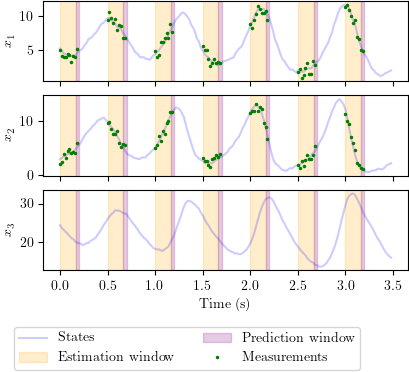}
	\end{center}
	\vspace{-0.5cm}
	\caption{Partially measured Lorenz attractor system with process and measurement noise,
 for $\theta = \begin{bmatrix}10 & 30 \end{bmatrix}^\top$.
	}\label{fig:LorenzMeas}
\end{figure}

\begin{figure}
	\vspace{-0.1cm}
	\begin{center}
		\includegraphics[width=0.95\linewidth]{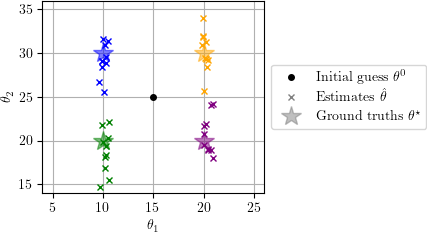}
	\end{center}
	\vspace{-0.5cm}
	\caption{Comparison between ground truth model parameters and estimates based on simulation data of the Lorenz attractor.
	}\label{fig:LorenzEstimates}
 \vspace{-0.5cm}
\end{figure}

As a study case of a nonlinear system, we consider the Lorenz attractor,
a system that can exhibit bounded chaotic behavior,
introduced by Edward Lorenz in 1963.
It has been used many times as a parameter estimation study case~\cite{Bonilla2008, annan2004efficient}.
The dynamics of the Lorenz attractor is governed by a set of differential equations of the form:
\begin{equation}\label{equation_Lorenz}
	\begin{aligned}
		\dot{x}_1(t) &= \theta_1 (x_2 - x_1), \\
		\dot{x}_2(t) &= x_1 (\theta_2 - x_3) - x_2, \\
		\dot{x}_3(t) &= x_1 x_2 - 2  x_3,
	\end{aligned}
\end{equation}
where $\theta \in \rr^2$ is the parameter vector and $x$ represents the state variables.
A trajectory generated from this system is depicted in Fig.~\ref{fig_Lorenz},
where the state switches between two orbits.
We consider discrete measurements of the two first states of the system,
taken at a sampling time of $\Delta t = 0.02$.
These measurements are subject to uniformly distributed noise in $[-1 \; 1]$, which is denoted by
$v_k \sim \Uu(-1, 1)$.
Additionally, we introduce time-discrete process noise at each time step, also uniformly distributed as $w_k \sim \mathcal{U}\left(-\frac{1}{4}, \frac{1}{4}\right)$.
We employ the Runge-Kutta-4 integration scheme for numerical discretization.
Finally, as described in Section~\ref{section_problem},
the available measurements consist of numerous subsequences, each with a length of $m=10$.
These subsequences, denoted as $y_{1, k}^{(i)}$ and $y_{2, k}^{(i)}$, are obtained from evenly spaced time steps $t^{(i)}$ with a step size $t_s = 0.5 \textup{sec.}$,
which is summarized by the following equations
\begin{equation}\label{equation_LorenzMeas}
	\begin{aligned}
		y_{1, k}^{(i)} &= x_1(i \cdot t_s + k\cdot \Delta t) + v_{1,k}^{(i)}, 
  &&\textup{for } k=1, \dots, m,
  \\
        y_{2, k}^{(i)} &= x_2(i \cdot t_s + k\cdot \Delta t) + v_{2,k}^{(i)},
		&&\textup{for } k=1, \dots, m,
	\end{aligned}
\end{equation}
where $v_{1,k}^{(i)}, v_{2,k}^{(i)} \sim \Uu(-1, 1)$.
We consider four different ground truth parameter vectors, $\theta^\star$, and generate $50$ measurement trajectories for each of them, with a random initial state $x(0) \sim \mathcal{U}(-10, 10)$. Each trajectory is of length $T=3.5 \textup{sec.}$, and Fig. \ref{fig:LorenzMeas} displays one such trajectory. Collating data from the 50 trajectories yields a dataset comprising $N=350$ subsequences for each $\theta^\star$.
For each dataset, we apply our parameter estimation approach and compare the estimates $\hat{\theta}$ with their corresponding ground truth $\theta^\star$.
To assess the robustness of our method, we repeat this procedure ten times with different noise realizations.
The initial guess for the optimization solver is set at $\theta^0 = {\begin{pmatrix} 15 & 25 \end{pmatrix}}\!^\top\!$.
The results are presented in Fig. \ref{fig:LorenzEstimates}.

These results demonstrate that across ten experimental realizations, the four ground-truth parameters are consistently recovered, rather accurately. In conclusion, this example showcases the capability of our proposed approach to accurately recover model parameters, even when based solely on noisy partial measurements of a chaotic system.

%% file: sections/conclusion.tex
The present paper introduces a novel extension of the prediction error method framework
to nonlinear systems, employing the moving horizon estimation method.
A theoretical property is established for the special case of linear time-invariant systems, and
numerical examples demonstrate the method's applicability to both linear and nonlinear systems.
These examples underscore the versatility of the proposed method and its performance in terms of recovering true system parameters.

The implementation of optimization algorithms with parallel computing 
is made possible by the design of the method.
Thus, it is an important topic for future research.
Another significant milestone for future research is
a comprehensive proof of asymptotic consistency under suitable conditions since it remains only partially proven in the present work.
Additionally, exploring variations in the method formulation looks promising.
For example, minimizing the multi-step-ahead prediction errors may increase the capabilities of the method.
As a conclusion, through rigorous theoretical analysis and empirical validation,
our work marks a solid starting point
for the topic of tuning a moving horizon state estimator for the identification of a nonlinear system,
which offers a versatile and reliable framework for parametric nonlinear system identification.

%% file: Simpson2024ParallelizableArXiv.bbl
\begin{thebibliography}{10}
\providecommand{\url}[1]{#1}
\csname url@samestyle\endcsname
\providecommand{\newblock}{\relax}
\providecommand{\bibinfo}[2]{#2}
\providecommand{\BIBentrySTDinterwordspacing}{\spaceskip=0pt\relax}
\providecommand{\BIBentryALTinterwordstretchfactor}{4}
\providecommand{\BIBentryALTinterwordspacing}{\spaceskip=\fontdimen2\font plus
\BIBentryALTinterwordstretchfactor\fontdimen3\font minus
  \fontdimen4\font\relax}
\providecommand{\BIBforeignlanguage}[2]{{%
\expandafter\ifx\csname l@#1\endcsname\relax
\typeout{** WARNING: IEEEtran.bst: No hyphenation pattern has been}%
\typeout{** loaded for the language `#1'. Using the pattern for}%
\typeout{** the default language instead.}%
\else
\language=\csname l@#1\endcsname
\fi
#2}}
\providecommand{\BIBdecl}{\relax}
\BIBdecl

\bibitem{Rawlings2017}
J.~B. Rawlings, D.~Q. Mayne, and M.~M. Diehl, \emph{Model Predictive Control:
  Theory, Computation, and Design}, 2nd~ed.\hskip 1em plus 0.5em minus
  0.4em\relax Nob Hill, 2017.

\bibitem{Verhaegen1994}
M.~Verhaegen, ``Identification of the deterministic part of {MIMO} state space
  models given in innovations form from input-output data,'' \emph{Automatica},
  vol.~30, pp. 61--74, 1994.

\bibitem{Van1994}
P.~Van~Overschee and B.~De~Moor, ``{N4SID}: Subspace algorithms for the
  identification of combined deterministic-stochastic systems,''
  \emph{Automatica}, vol.~30, pp. 75--93, 1994.

\bibitem{Ljung1999}
L.~Ljung, \emph{{S}ystem identification: Theory for the User}.\hskip 1em plus
  0.5em minus 0.4em\relax Upper Saddle River, N.J.: Prentice Hall, 1999.

\bibitem{Ljung2002}
------, ``Prediction error estimation methods,'' \emph{Circuits, Syst. and
  Signal Process.}, vol.~21, pp. 11--21, 2002.

\bibitem{Simpson2023}
L.~Simpson, A.~Ghezzi, J.~Asprion, and M.~Diehl, ``An efficient method for the
  joint estimation of system parameters and noise covariances for linear
  time-variant systems,'' in \emph{Proc. IEEE Conf. Decision and Control
  (CDC)}, 2023, pp. 4524--4529.

\bibitem{Kalman1960}
R.~Kalman, ``A new approach to linear filtering and prediction problems,''
  \emph{Trans. of the ASME--J. of Basic Eng.}, vol.~82, pp. 35--45, 1960.

\bibitem{McGee1985}
L.~A. McGee and S.~F. Schmidt, ``Discovery of the {K}alman filter as a
  practical tool for aerospace and industry,'' \emph{NASA Technical
  Memorandum}, p.~21, 1985.

\bibitem{alessandri2008moving}
A.~Alessandri, M.~Baglietto, and G.~Battistelli, ``Moving-horizon state
  estimation for nonlinear discrete-time systems: New stability results and
  approximation schemes,'' \emph{Automatica}, vol.~44, pp. 1753--1765, 2008.

\bibitem{Baumgaertner2019}
K.~Baumg{\"a}rtner, A.~Zanelli, and M.~Diehl, ``Zero-order moving horizon
  estimation,'' in \emph{Proc. IEEE Conf. Decision and Control (CDC)}, 2019.

\bibitem{schiller2023lyapunov}
J.~D. Schiller, S.~Muntwiler, J.~K{\"o}hler, M.~N. Zeilinger, and M.~A.
  M{\"u}ller, ``A {Lyapunov} function for robust stability of moving horizon
  estimation,'' \emph{IEEE Trans. Automatic Control}, 2023.

\bibitem{zheng2020}
Y.~Zheng and N.~Li, ``Non-asymptotic identification of linear dynamical systems
  using multiple trajectories,'' \emph{IEEE Contr. Sys. Lett.}, vol.~5, pp.
  1693--1698, 2020.

\bibitem{Valluru2017}
J.~Valluru, P.~Lakhmani, S.~C. Patwardhan, and L.~T. Biegler, ``Development of
  moving window state and parameter estimators under maximum likelihood and
  {B}ayesian frameworks,'' \emph{J.\ Proc.\ Contr.}, vol.~60, pp. 48--67, 2017.

\bibitem{astrom1979maximum}
K.~Astrom, ``Maximum likelihood and prediction error methods,'' \emph{Proc.
  IFAC World Congress}, vol.~12, pp. 551--574, 1979.

\bibitem{tu2022learning}
S.~Tu, R.~Frostig, and M.~Soltanolkotabi, ``Learning from many trajectories,''
  \emph{arXiv preprint arXiv:2203.17193}, 2022.

\bibitem{dean2020sample}
S.~Dean, H.~Mania, N.~Matni, B.~Recht, and S.~Tu, ``On the sample complexity of
  the linear quadratic regulator,'' \emph{Found. of Comput. Math.}, vol.~20,
  pp. 633--679, 2020.

\bibitem{Muntwiler2022}
S.~Muntwiler, K.~P. Wabersich, and M.~N. Zeilinger, ``Learning-based moving
  horizon estimation through differentiable convex optimization layers,'' in
  \emph{Learning for Dynamics and Control Conference}.\hskip 1em plus 0.5em
  minus 0.4em\relax PMLR, 2022, pp. 153--165.

\bibitem{esfahani2021reinforcement}
H.~N. Esfahani, A.~B. Kordabad, and S.~Gros, ``Reinforcement learning based on
  {MPC}/{MHE} for unmodeled and partially observable dynamics,'' in \emph{Proc.
  Amer. Control Conf. (ACC)}.\hskip 1em plus 0.5em minus 0.4em\relax IEEE,
  2021, pp. 2121--2126.

\bibitem{esfahani2023learning}
H.~N. Esfahani, A.~B. Kordabad, W.~Cai, and S.~Gros, ``Learning-based state
  estimation and control using {MHE} and {MPC} schemes with imperfect models,''
  \emph{Europ.\ J.\ Contr.}, p. 100880, 2023.

\bibitem{rao2002constrained}
C.~V. Rao and J.~B. Rawlings, ``Constrained process monitoring: Moving-horizon
  approach,'' \emph{AIChE J.}, vol.~48, pp. 97--109, 2002.

\bibitem{Waechter2006}
A.~W\"achter and L.~T. Biegler, ``On the implementation of an interior-point
  filter line-search algorithm for large-scale nonlinear programming,''
  \emph{Math. Program.}, vol. 106, pp. 25--57, 2006.

\bibitem{Andersson2019}
J.~A.~E. Andersson, J.~Gillis, G.~Horn, J.~B. Rawlings, and M.~Diehl,
  ``{CasADi} -- a software framework for nonlinear optimization and optimal
  control,'' \emph{Math. Program. Comput.}, vol.~11, pp. 1--36, 2019.

\bibitem{Bonilla2008}
J.~Bonilla, M.~Diehl, and B.~D. Moor, ``A nonlinear least squares estimation
  procedure without initial parameter guesses,'' \emph{Proc. IEEE Conf.
  Decision and Control (CDC)}, 2008.

\bibitem{annan2004efficient}
J.~Annan and J.~Hargreaves, ``Efficient parameter estimation for a highly
  chaotic system,'' \emph{Tellus}, vol.~56, pp. 520--526, 2004.

\end{thebibliography}
